\documentclass{article}
\usepackage{amssymb}

\usepackage{graphicx}
\usepackage{amsmath}
\usepackage{color}
\usepackage{float}
\usepackage{array}

\newtheorem{theorem}{Theorem}[section]

\newtheorem{corollary}[theorem]{Corollary}

\newtheorem{lemma}[theorem]{Lemma}

\newtheorem{problem}[theorem]{Problem}
\newtheorem{proposition}[theorem]{Proposition}
\newtheorem{remark}[theorem]{Remark}

\newenvironment{proof}[1][Proof]{\textbf{#1.} }{\ \rule{0.5em}{0.5em}}
\def\b{\beta}

\def\O{\Omega}

\def\eps{\varepsilon}

\def\Om{\Omega}

\newcommand{\R}{\mathbb{R}}

\newcommand{\cqfd}

\newcommand{\jj}{j_{1,1}}

\newcommand{\dn}{\partial_\nu}

\newcommand{\bOm}{\partial\Omega}

\def\div{\mbox{div}}

\DeclareMathOperator{\dist}{dist}

\title{ A Talenti-type  result for the Pompeiu-Schiffer Conjecture  }
\author{ Diaraf SECK$^{1,\,}$ \footnote{diaraf.seck@ucad.edu.sn}\\\\
$^{1}$ Laboratoire de Math\'ematiques de la D\'ecision et \\
d'Analyse Num\'erique (LMDAN), BP 16889, Dakar Fann, S\'enegal\\
Ecole Doctorale de Math\'ematiques et Informatique U.C.A.D. Dakar,  S\'en\'egal.}

\begin{document}
\maketitle
\begin{abstract}
Schiffer's conjecture asserts that the ball is the only bounded simply connected domain carrying a
Neumann eigenfunction of the Laplacian that is constant on the boundary; it is equivalent to the
Pompeiu problem. In August 2026 it was disproved in the plane by computer-assisted constructions.
We revisit the shape-optimization approach of~\cite{Se24} in the light of these counterexamples.
We first show that the vanishing of the material derivative on the boundary, assumed in~\cite{Se24},
never holds, and we compute the correct optimality condition, involving an adjoint state, for the
associated functional $J_2$. We prove that every planar Schiffer domain satisfies $\max|u| > \sqrt{2}|c|.$ Our main result is
a Talenti-type theorem: if $u < c$ in $\Omega$ , then $kR^* \geq  j_{1,1},$  where $R^*$ is the radius of the disc of
equal area, and equality holds only for the disc. Finally, using Siegel's transcendence
theorem, we prove the infinitesimal rigidity of the disc for the Berenstein problem.
\end{abstract}
{\bf Keywords:}  Schiffer conjecture, Pompeiu problem, 
overdetermined problem, Talenti rearrangement, $P$-functions,  rigidity, Bessel
functions,  simple eigenvalue, symmetry.\\
\section{Introduction} 

Let $\Omega \subset\R^N$ be a bounded, regular, simply connected domain. Consider the overdetermined problem
\begin{equation}\label{eq:schiffer}
\Delta u+\lambda u=0 \ \text{in }\Omega,\qquad \partial_{\nu} u=0\ \text{and}\ u=c\ \text{on }\partial \Omega,
\qquad \lambda>0,\ u\not\equiv0 .
\end{equation}
Such a domain satisfying (\ref{eq:schiffer}) is called \emph{Schiffer's domain.}\\
\emph{Schiffer's conjecture} asserted that the existence of a solution of~\eqref{eq:schiffer} forces $\Omega$
to be a ball and $u$ to be radial. The \emph{Pompeiu problem}~\cite{Po29} asks whether a continuous
function $f$ is determined by its integrals over all rigid motions of $\Om$; the ball fails this property.
Brown, Schreiber and Taylor~\cite{BST73} and Williams~\cite{Wi76} showed that, for smooth contractible
domains, failure of the Pompeiu property is equivalent to the existence of a solution
of~\eqref{eq:schiffer}: Green's formula then shows that the Fourier transform of $\mathbf 1_\Omega$
vanishes on the sphere $|\xi|=\sqrt\lambda$. This is Problem~80 in Yau's list~\cite{Ya82}. We call $\Omega$ a
\emph{Schiffer domain} if~\eqref{eq:schiffer} has a solution.

\subsubsection*{Partial results} Berenstein and Yang~\cite{BY87} proved that infinitely many admissible
eigenvalues force the ball; Williams~\cite{Wi81} that the boundary of a Schiffer domain is analytic. The
conjecture is known for small eigenvalues~\cite{Av86,De12}, for quadrature domains~\cite{Eb93}, for some
centrally symmetric convex domains~\cite{Mo26,DSWZ25}, when the third normal derivative is
constant~\cite{Li07}, and the regular pentagon is not a Schiffer domain~\cite{NSY20}. Counterexamples were
known in other geometries~\cite{FMW25,CLF25}.

\subsubsection*{The 2026 disproof} Colbrook and Stepaniants~\cite{CS26} constructed, through a
computer-assisted proof, a bounded simply connected planar domain with analytic boundary and $D_{10}$
symmetry, which is not a disc and carries a solution of~\eqref{eq:schiffer} with
$k=\sqrt\lambda\in(31.967007261,\,31.967007293)$; an independent construction is due to Cao-Labora and
de Dios Pont~\cite{CLDP26}. Convex counterexamples in dimensions two to eighteen have been
announced~\cite{Cvx26}. In the same spirit, Colbrook, Sadeghi and Stepaniants~\cite{CSS26} disproved the
\emph{Berenstein conjecture}: there is a non-circular planar domain ($D_{13}$ symmetry) carrying a solution of
\begin{equation}\label{eq:berenstein}
\Delta u+k^2u=0\ \text{in }\O,\qquad u=0\ \text{and}\ \partial_{\nu} u=1\ \text{on }\partial \Omega .
\end{equation}
\subsection{The shape-optimization approach}
The link between these conjectures and shape optimization goes back in particular to~\cite{Gr00}, where
the ball is obtained as a stable minimizer of a shape functional, and to~\cite{CH99}, where Schiffer
domains appear as critical points of an eigenvalue under a volume constraint. In~\cite{Se24}, the author
distinguishes Conjectures~1--4 (Pompeiu, Schiffer, spherical zeros of $\widehat{\mathbf 1_\Om}$, all
equivalent) from Conjecture~5, which is problem~\eqref{eq:berenstein}. Three kinds of results are obtained
there: a Serrin-type symmetry theorem~\cite{Se71} on a domain with a hole, assuming that $-\Delta-k^2$
satisfies the maximum principle; Conjecture~5 in the regime $|k|<j_{n,1}/R$ (Theorem~4.2 of~\cite{Se24}), by
shape minimization and monotonicity of the Lagrange multiplier; and a Riemannian approach~\cite{MM06}
to the Pompeiu problem through a functional
\[
J_2(\Omega)=\int_{\partial\Omega}u^2\,d\sigma+\gamma|\Omega|
\]
constrained by the Neumann problem, whose analysis relies on the assumption that the material derivative
of $u$ vanishes on the boundary. The limits of the method are identified there: the maximum principle only
holds for $k$ small compared with the radius.
\subsection{Results}
The present work re-examines~\cite{Se24} after 2026.
And we  have the following   results: 
\begin{itemize}
\item[(i)] the vanishing of the material derivative on the boundary, assumed in \cite{Se24}, never holds; we compute the correct optimality condition, with an adjoint state, for the boundary functional of \cite{Se24};
\item[(ii)] a Talenti-type theorem: if $u<c$ in
$\Omega$, then $kR^*\geq j_{1, 1}.$ If we have equality then  $\Omega$ is a disc; 
\item[(iii)] bounds valid at every
frequency ($\max|u|>\sqrt2\,|c|$, $|\nabla u|^2+\lambda u^2\le\lambda\max u^2$). We also prove the
infinitesimal rigidity of the disc. 
\end{itemize}

The common thread of these results is the condition of
sign  (H1) $  u < c$ in $\Omega$: as in~\cite{Se24}, where the maximum
principle requires $|k|<j_{n,1}/R$, rigidity can be proved as long as a radial test function keeps a constant  sign, whereas
the counterexamples live at high frequency.
\section{Preliminaries}\label{sec:prelim}

Let $\Om$ be a Schiffer domain and $u$ a solution of~\eqref{eq:schiffer} with $c\neq0$ (the case $c=0$ is
excluded by uniqueness for the Cauchy problem). Since $u=c$ and $\dn u=0$ on $\bOm$,
\begin{equation}\label{eq:gradzero}
\nabla u=0\quad\text{on }\bOm .
\end{equation}


Let us introduce   the oriented distance, called also signed distance defined as follows, for $x \in \R^2:$
 \begin{equation*}
 s (x)= 
    \left\{\begin{array}{ccc}
  -\dist(x,\partial\Omega) & if & x\in\Omega,\\ 
    0 &if & x \in \partial \Omega,\\
 \dist(x,\partial\Omega) &if & x\notin \Omega
                  \end{array}
\right.
\end{equation*}
We denote by $\nu$ the outward normal and by $H=\div \nu = \Delta s$  the mean curvature (see \cite{DeZo}); in dimension~2 we write $\kappa$ for the curvature.


Thus  if $s < 0$ one is  in $\Omega$, for  $s = 0 $ the position is  on $\partial \Omega$ and  when $s> 0$, we are outside of $\Omega$. 
 If $\partial \Omega$ is of class $\mathcal C^2$, there exists  $\delta >0$  (in dimension 2 one can take  $\delta <1/\max |\kappa|,$ provided  $\partial \Omega$ does not come close to itself) such that every $x$ with $|s(x)|< \delta$  can be written uniquely as

\begin{eqnarray*}
x = y + s\,\nu(y), \qquad y\in\partial\Omega \ \text{(orthogonal projection of } x\text{)},\quad s = s(x).
\end{eqnarray*}

The pair $(y, s)$ forms the \textbf{normal coordinates}. The level sets $ \Gamma_{s} = \{s = const\}$  are the curves parallel to $\partial \Omega = \Gamma_0$.
\begin{figure}[htbp]\centering
\includegraphics[width=0.62\linewidth]{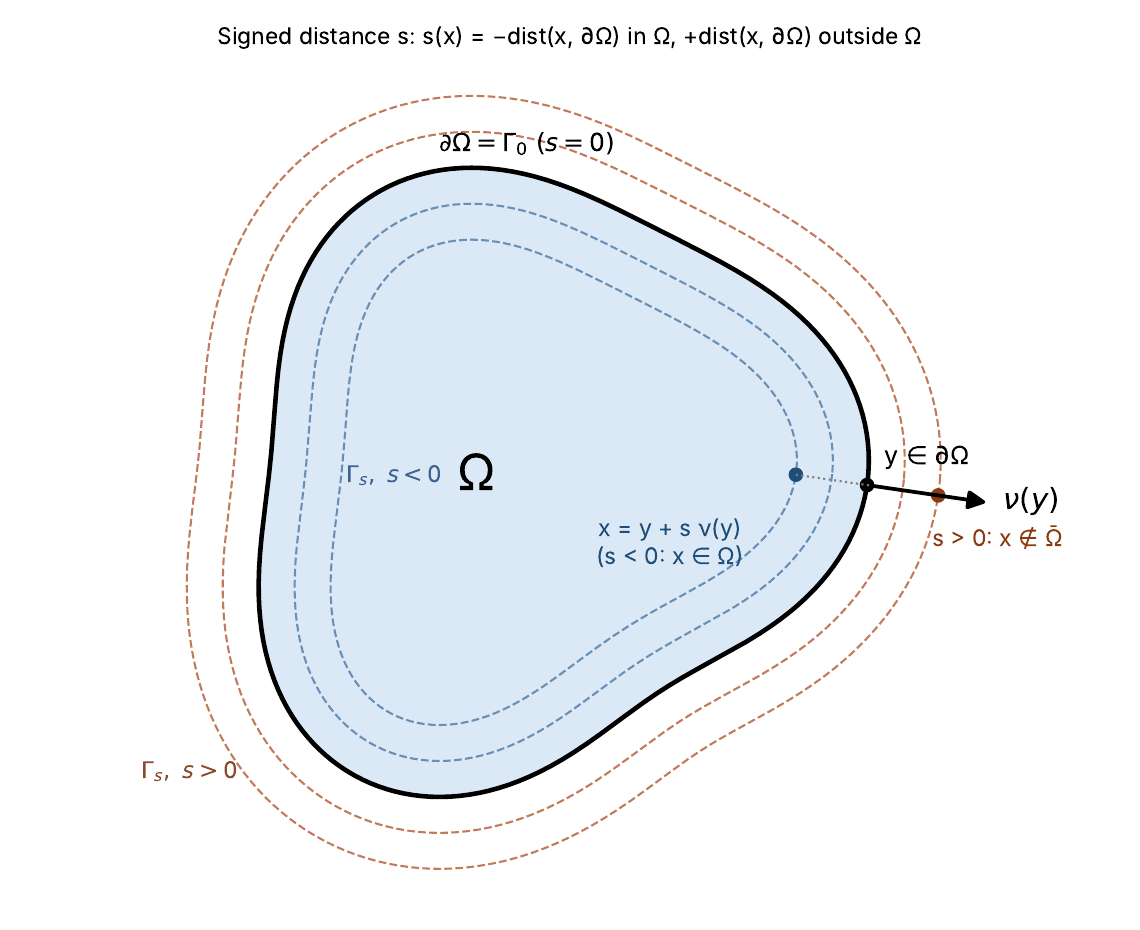}
\caption{Signed distance s, normal coordinates (y, s) and parallel curves $\Gamma_s: s< 0$ in $\Omega$ , $s >0$ outside} \label{fig:s}
\end{figure}
\newpage
Let us point out that  in this neighborhood, $s$ is smooth and:

\begin{enumerate}
\item  $|\nabla s| = 1,\qquad \nabla s(x) = \nu(y),\qquad \partial_s = \nu\cdot\nabla,\qquad \Delta s = \operatorname{div}\nu = H_s.$
Here $H_s$ is the mean curvature of $\Gamma_s$ ; on $ \partial \Omega , H = \kappa$ in dimension 2.
\item  In dimension 2, $\kappa_s = \kappa /(1 + s\kappa),$ hence $\partial_{s} \kappa_s =- \kappa ^2$   at $s = 0.$  For the disc of radius $R, s = r -R$ and $\Delta s = 1/r = \kappa/(1 + s\kappa).$
\end{enumerate}

\begin{lemma}[Normal derivatives]\label{lem:normal}
On $\bOm$,
\[
\dn^2u=-\lambda c,\qquad \dn^3u=\lambda c\,H .
\]
In dimension~2, moreover, $\dn^4u=\lambda c\,(\lambda-3\kappa^2)$.
\end{lemma}

\begin{proof}
Near $\bOm$, $\Delta u=\partial_s^2u+H_s\,\partial_su+\Delta_{\Gamma_s}u$, where $\Gamma_s$ is the parallel
surface at distance $s$ and $H_s$ its mean curvature. On $\bOm$, $\partial_su=0$ and $\Delta_\Gamma u=0$
since $u$ is constant; the equation gives $\partial_s^2u=-\lambda c$. Since the first three coefficients of the
expansion in $s$ do not depend on the boundary point, the term $\Delta_{\Gamma_s}u$ is $O(s^3)$;
differentiating the equation in $s$ and setting $s=0$ gives $\partial_s^3u+H\partial_s^2u=0$. In
dimension~2, parallel curves satisfy $\partial_s\kappa_s=-\kappa^2$, and one more differentiation gives the
last formula.
\end{proof}

\begin{corollary}[Liu ~\cite{Li07}]
$\dn^3u$ is constant on $\bOm$ if and only if $\Om$ is a ball.
\end{corollary}
\begin{proof}
By the lemma this is equivalent to $H$ being constant; apply Alexandrov's theorem.
\end{proof}

\begin{remark}
These local relations do not constrain $\Om$: since the boundary is analytic~\cite{Wi81}, the
Cauchy--Kovalevskaya theorem provides, for every analytic curve, a local solution near the boundary with
Cauchy data $(c,0)$. The obstruction is global: this solution must extend to the whole of $\Om$.
\end{remark}

\begin{lemma}[Pohozaev-Integral identities]\label{lem:rellich}
For a Schiffer domain in $\R^N$,
\[
\int_\Om u=0,\qquad \int_\Om|\nabla u|^2=\lambda\int_\Om u^2,\qquad \int_\Om u^2=\frac N2\,c^2|\Om| .
\]
\end{lemma}
\begin{proof}
The first identity follows by integrating the equation, the second by multiplying it by $u$. For the
third, multiply the equation by $x\cdot\nabla u$ (Rellich--Pohozaev identity): by~\eqref{eq:gradzero},
every boundary term containing $\nabla u$ vanishes, and what remains is
$\lambda\big(\tfrac N2\int_\Om u^2-\tfrac12\int_{\bOm}u^2\,x\cdot\nu\big)=\big(\tfrac N2-1\big)\int_\Om|\nabla u|^2$,
which gives the claim since $\int_{\bOm}x\cdot\nu=N|\Om|$.
\end{proof}

These identities hold for every Schiffer domain, counterexamples included; they therefore do not
characterize the ball.

\section{Material derivative and the functional $J_2$}\label{sec:material}

Let $V$ be a deformation field, $\Om_t=(I+tV)(\Om)$, $\alpha=V\cdot\nu$ on $\bOm$, and $u_t$ the normalized
Neumann eigenfunction of $\Om_t$ (simple eigenvalue). We denote by $\dot u$ the material derivative and by
$u'$ the shape derivative, so that $\dot u=u'+\nabla u\cdot V$.

For the Dirichlet problem, the condition $u_t=0$ on $\bOm_t$ is transported by the deformation, hence
$\dot u=0$ on the boundary. For the Neumann problem, it is $\partial_{\nu_t}u_t=0$ that is transported, not
the value of $u_t$.

\begin{proposition}\label{prop:material}
Let $\Om$ be a Schiffer domain (simple eigenvalue, $\|u\|_{L^2}=1$). It is not true that $\dot u=0$ on
$\bOm$ for every deformation preserving the volume to first order.
\end{proposition}

\begin{proof}
By~\eqref{eq:gradzero}, $\dot u=u'$ on $\bOm$. The shape derivative satisfies
\[
\dn u'=\div_\Gamma(\alpha\nabla_\Gamma u)+\lambda u\,\alpha=\lambda c\,\alpha,\qquad
\lambda'=\int_{\bOm}\big(|\nabla_\Gamma u|^2-\lambda u^2\big)\alpha=-\lambda c^2\int_{\bOm}\alpha .
\]
Let $w$ solve $\Delta w+\lambda w=0$ in $\Om$. Since $-\Delta u'-\lambda u'=\lambda'u$, Green's formula gives
$-\lambda'\int_\Om uw=\int_{\bOm}(\dn u'\,w-u'\,\dn w)$. If $u'=0$ on $\bOm$ and $\int\alpha=0$, then
$\lambda'=0$ and $\int_{\bOm}\alpha\,w=0$ for every $\alpha$ with zero mean, so $w$ is constant on $\bOm$.
Taking $w=e^{ik\omega\cdot x}$ with $k^2=\lambda$, the function $\omega\cdot x$ would take values in a
discrete set on the connected boundary, hence would be constant on $\bOm$ for every direction $\omega$,
which is absurd.
\end{proof}

In particular, the derivative of $J_2$ at a Schiffer domain reads
$dJ_2(\Om)[V]=\int_{\bOm}\big(2c\,u'+\alpha(Hc^2+\gamma)\big)$, and the term $2c\,u'$ cannot be dropped.
If it were, criticality would give constant $H$, hence a ball: the assumption $\dot u=0$ already contains
the conclusion.

\begin{proposition}[Optimality condition with adjoint state]\label{prop:adjoint}
Let $u$ be the normalized eigenfunction ($\int_\Om u^2=1$), $\lambda$ simple, $\beta=\int_{\bOm}u^2$  ($\beta$ is chosen so that the Fredholm compatibility condition is
satisfied), and
let $p$ be the adjoint state defined by
\[
-\Delta p-\lambda p=-\beta u\ \text{in }\Om,\qquad \dn p=u\ \text{on }\bOm,\qquad \int_\Om pu=0 .
\]
Then $dJ_2(\Om)[\alpha]=\int_{\bOm}\alpha\,\big(Hu^2+\gamma-\beta u^2+2\lambda u\,p-2\nabla_\Gamma u\cdot\nabla_\Gamma p\big)$.
If $\Om$ is a planar Schiffer domain, $p=-\frac{\beta}{2\lambda}\,x\cdot\nabla u+c\,q$, where $q$  (defined up to the addition of a multiple of u) solves
\[
\Delta q+\lambda q=0\ \text{in }\Om,\qquad \dn q=1-\frac{|\bOm|}{2|\Om|}\,x\cdot\nu\ \text{on }\bOm ,
\]
and $\Om$ is a critical point of $J_2$ (for some $\gamma$) if and only if
\[
T:=H+2\lambda q\ \text{is constant on }\bOm .
\]
\end{proposition}

\begin{proof}[Sketch of proof]
The general formula follows from Hadamard's formula and Green's formula applied to the pair $(u',p)$. At
a Schiffer domain one uses $(-\Delta-\lambda)(x\cdot\nabla u)=2\lambda u$ and $D^2u=-\lambda c\,\nu\otimes\nu$
on $\bOm$, which gives $\dn(x\cdot\nabla u)=-\lambda c\,x\cdot\nu$; the identities $\beta=c^2|\bOm|$ and
$c^2|\Om|=1$ (Lemma~\ref{lem:rellich}) give the Neumann condition for $q$. Finally $u=c$, $\nabla_\Gamma u=0$
and $p=c\,q$ on $\bOm$.
\end{proof}
\section{The function $P=|\nabla u|^2+\lambda u^2$ in dimension 2}\label{sec:P}

Throughout this section, $\Om\subset\R^2$ is a Schiffer domain and $P=|\nabla u|^2+\lambda u^2$.

\begin{lemma}[Exact identity]\label{lem:Pexact}
On the open set $\{\nabla u\neq0\}$,
\begin{equation}\label{eq:Pexact}
\Delta P+\frac{2\lambda u}{|\nabla u|^2}\,\nabla u\cdot\nabla P-\frac{|\nabla P|^2}{|\nabla u|^2}=0 .
\end{equation}
\end{lemma}
This identity is due to the $P-$ functions theory developped by Payne, Philippin and Sperb. For more details see for instance \cite{Sp81}.

\begin{proof}
We have $\nabla P=2D^2u\,\nabla u+2\lambda u\nabla u$ and, by Bochner's formula and the equation,
$\Delta P=2|D^2u|^2-2\lambda^2u^2$. Let $\rho=|\nabla u|$, $e_1=\nabla u/\rho$, $e_2\perp e_1$, and write
$D^2u=\big(\begin{matrix}a&b\\b&d\end{matrix}\big)$ in this frame, with $a+d=-\lambda u$. Setting
$A=a+\lambda u$ we get $d=-A$ (this is where dimension~2 is used), then $\nabla P=2\rho(Ae_1+be_2)$, hence
$\nabla u\cdot\nabla P=2\rho^2A$ and $|\nabla P|^2=4\rho^2(A^2+b^2)$. Moreover
$|D^2u|^2=(A-\lambda u)^2+2b^2+A^2$, so
$\Delta P=4A^2+4b^2-4\lambda uA=|\nabla P|^2/\rho^2-2\lambda u\,\nabla u\cdot\nabla P/\rho^2$.
\end{proof}

\begin{remark}
In dimension $N\ge3$, the orthogonal complement of $\nabla u$ contains an $(N-1)\times(N-1)$ block and one
only gets an inequality, as in Payne's theory of $P$-functions~\cite{Sp81}.
\end{remark}

Identity~\eqref{eq:Pexact} reads $\Delta P+B\cdot\nabla P=0$ with $B=(2\lambda u\nabla u-\nabla P)/|\nabla u|^2$,
locally bounded on each component of $\{\nabla u\neq0\}$. Before proceeding further let us point out that $\partial \Omega$ is part of the critical set of $u.$ This gives the following principle.

\begin{corollary}\label{cor:extremaP}
The extrema of $P$ on $\overline\Om$ are attained on $\bOm$, where $P=\lambda c^2$, or at critical points of
$u$, where $P=\lambda u^2$.
\end{corollary}

\begin{lemma}[$P$ near the boundary]\label{lem:Pboundary}
With $s<0$ inside, $P=\lambda c^2-\frac23\lambda^2c^2\,\kappa\,s^3+O(s^4)$. Just inside, $P-\lambda c^2$ has
the sign of $\kappa$.
\end{lemma}
\begin{proof}
This follows from Lemma~\ref{lem:normal} and $|\nabla_\Gamma u|^2=O(s^6)$.
\end{proof}

\begin{proposition}\label{prop:cp}
Let $\Om\subset\R^2$ be a Schiffer domain. 
\begin{enumerate}\label{roman*}
\item There is a critical point $x_0$ of $u$ with $|u(x_0)|>|c|$.
\item If $\Om$ is not convex, there is a critical point $x_1$ of $u$ with $|u(x_1)|<|c|$.
\end{enumerate}
In particular, if $u$ has a single critical point in $\Om$, then $\Om$ is convex (in fact a disc, by \cite{WG94}).
\end{proposition}

\begin{proof}
(i) Since $\int_{\bOm}\kappa=2\pi$, there is a boundary point where $\kappa>0$; by
Lemma~\ref{lem:Pboundary}, $P>\lambda c^2$ near it, so $M=\max P$ is attained inside. If the maximum point
is not critical, the strong maximum principle gives $P\equiv M$ on the component $W$ of $\{\nabla u\neq0\}$
containing it; $W$ does not touch $\bOm$ (where $P=\lambda c^2<M$), so its boundary consists of critical
points, where $\lambda u^2=P=M>\lambda c^2$. (ii) Same argument with the minimum, at a boundary point
where $\kappa<0$. Finally, a single critical point is the extremum of $|u|$; it satisfies (i), which rules out (ii).
\end{proof}

\begin{theorem}\label{TheoremC}
Let $\Om\subset\R^2$ be a Schiffer domain and $M=\max_{\overline\Om}|u|$. Then
\begin{enumerate}\label{alph*}
\item $M>\sqrt2\,|c|$;
\item $|\nabla u|^2+\lambda u^2\le\lambda M^2$ in $\Om$;
\item if $|u(x^*)|=M$, then $\dist(x^*,\bOm)\ge\big(\frac\pi2\mp\arcsin\frac{|c|}{M}\big)/k$, with the sign $-$
if $u(x^*)$ and $c$ have the same sign and $+$ otherwise.
\end{enumerate}
\end{theorem}

\begin{proof}
(b) By Corollary~\ref{cor:extremaP}, $\max P\le\max(\lambda c^2,\lambda M^2)=\lambda M^2$, since $M>|c|$ by
Proposition~\ref{prop:cp}. (a) By Lemma~\ref{lem:rellich} in dimension~2,
$\int_\Om P=2\lambda\int_\Om u^2=2\lambda c^2|\Om|$: the mean of $P$ equals $2\lambda c^2$ and is strictly less
than $\lambda M^2$ because $P$ is not constant. (c) By (b), $w=\arcsin(u/M)$ satisfies $|\nabla w|\le k$; it
equals $\pm\pi/2$ at $x^*$ and $\arcsin(c/M)$ on $\bOm$.
\end{proof}
\section{A Talenti-type rigidity theorem}\label{sec:talenti}

We write $k=\sqrt\lambda$, $R^*=(|\Om|/\pi)^{1/2}$ for the radius of the disc with the same area, and
$\jj\approx3.8317$ for the first positive zero of $J_1$.

\begin{theorem}\label{TheoremB}
Let $\Om\subset\R^2$ be a bounded domain with $C^2$ boundary and $u\in C^2(\overline\Om)$ a solution
of~\eqref{eq:schiffer} with $c>0$. Assume
\[
\text{(H1)}\qquad u<c\ \text{in }\Om .
\]
Then $kR^*\ge\jj$. And if  $kR^*= \jj$ then $\Om$ is a disc.
\end{theorem}

\begin{proof}
We normalize $c=1$.

\emph{Step 0.} We have $\int_\Om u=0$ (Lemma~\ref{lem:rellich}); $u$ is analytic and non-constant, so
$\{\nabla u=0\}$ and each level set are negligible; finally $m_0:=\min u<0$.

\emph{Step 1 (role of (H1)).} For $t\in(m_0,1)$, let $E_t=\{u<t\}$. The set $\{u\le t\}$ is compact and does
not meet $\bOm$, so $E_t\Subset\Om$ and $\partial E_t=\{u=t\}$. For almost every $t$ (Sard), $\{u=t\}$ is a
finite union of smooth closed curves with outward normal $\nabla u/|\nabla u|$. Moreover
$\bigcup_{t<1}E_t=\Om$.

\emph{Step 2 (flux).} For regular $t$,
$F(t):=\int_{\{u=t\}}|\nabla u|=\int_{E_t}\Delta u=-\lambda\int_{E_t}u>0$.

\emph{Step 3 (distribution function).} By the coarea formula, $\mu(t)=|E_t|$ is absolutely continuous and
strictly increasing, with $\mu'(t)=\int_{\{u=t\}}|\nabla u|^{-1}>0$ almost everywhere, $\mu(m_0)=0$ and
$\mu(1^-)=|\Om|$.

\emph{Step 4 (Cauchy--Schwarz and isoperimetry).} If $L(t)$ is the length of $\{u=t\}$, then
$L(t)^2\le\mu'(t)F(t)$, and the planar isoperimetric inequality gives, for almost every $t$,
\begin{equation}\label{eq:talenti}
4\pi\mu(t)\le L(t)^2\le\mu'(t)F(t)=-\lambda\,\mu'(t)\int_{E_t}u .
\end{equation}
Without (H1), $\partial E_t$ would contain $\bOm$ for $t>1$ and~\eqref{eq:talenti} would fail.

\emph{Step 5 (rearrangement).} Let $\varphi:(0,|\Om|)\to(m_0,1)$ be the inverse of $\mu$; it is absolutely
continuous with $\varphi'(a)=1/\mu'(\varphi(a))$ almost everywhere. The function $u^*(x)=\varphi(\pi|x|^2)$ is the
increasing radial rearrangement of $u$. With $\Phi(a)=\int_0^a\varphi$ we have $\int_{E_t}u=\Phi(\mu(t))$,
and~\eqref{eq:talenti} becomes
\begin{equation}\label{eq:psi}
\psi(a):=4\pi a\,\varphi'(a)+\lambda\Phi(a)\le0\quad\text{for almost every }a\in(0,|\Om|),
\end{equation}
with $\varphi(|\Om|^-)=1$, $\Phi(0)=0$ and $\Phi(|\Om|)=\int_\Om u=0$. In the variable $a=\pi r^2$ the radial
Laplacian reads $\frac{d}{da}\big(4\pi a\frac{d}{da}\big)$: $\psi'=\Delta u^*+\lambda u^*$, and $\psi\equiv0$
expresses that $u^*$ solves the radial equation.

\emph{Step 6 (test function).} Let $Z(a)=J_0\big(k\sqrt{a/\pi}\big)$. With $a=\pi r^2$ and $(xJ_1)'=xJ_0$,
\[
Z'(a)=-\frac{kJ_1(kr)}{2\pi r},\qquad 4\pi aZ'(a)=-2\pi k\,rJ_1(kr),\qquad (4\pi aZ')'=-\lambda Z .
\]

\emph{Step 7 (exact identity).} Let $A=|\Om|$. Two integrations by parts ($\varphi$ bounded and absolutely
continuous, $aZ'(a)\to0$ at $0$) give
\begin{align*}
\int_0^A4\pi a\varphi'Z'&=4\pi AZ'(A)\varphi(A)+\lambda\int_0^A\varphi Z,\\
\lambda\int_0^A\Phi Z'&=\lambda\Phi(A)Z(A)-\lambda\int_0^A\varphi Z=-\lambda\int_0^A\varphi Z,
\end{align*}
hence
\begin{equation}\label{eq:I}
I:=\int_0^{A}\psi(a)\,Z'(a)\,da=4\pi AZ'(A)=-2\pi\,kR^*J_1(kR^*) .
\end{equation}

\emph{Step 8 (sign of $J_1$).} Assume $kR^*\le\jj$. For $0<a<A$ we have $kr<\jj$, so $Z'(a)<0$ and,
by~\eqref{eq:psi}, $I\ge0$. But $J_1(kR^*)\ge0$, so $I\le0$ by~\eqref{eq:I}. Thus $I=0$, which gives
$J_1(kR^*)=0$, that is $kR^*=\jj$, and $\psi=0$ almost everywhere.

\emph{Step 9 (equality).} If $\psi=0$ almost everywhere, the inequalities in~\eqref{eq:talenti} are
equalities for almost every $t$: $E_t$ is a disc. The sets $E_t$ form an increasing family of nested discs
whose radii tend to $R^*$ and whose centres converge ($|c_t-c_s|\le r_s-r_t$); their union $\Om$ is a disc.
Conversely, on the disc with $kR=\jj$, $u=J_0(kr)/J_0(\jj)$ is increasing in $r$ and satisfies (H1).
\end{proof}
\begin{problem}
 Does $(H_1)$ imply that $\Omega$ is a disc?
\end{problem}

\begin{remark}
The interval $(0,\jj)$ is exactly where $Z$ is monotone; beyond it, $Z'$ changes sign and the argument no
longer concludes. The positivity of $J_1$ plays here the role of the maximum principle in~\cite{Se24}. The
approach is that of Talenti's comparison principle~\cite{Ta76}; what is specific to the Schiffer problem is
the vanishing of the boundary terms in~\eqref{eq:I}, due to $\dn u=0$ and $\int_\Om u=0$. On Schiffer discs
of order $n\ge2$, (H1) fails ($\max u=3.33,\,1.61,\,4.58$ for $n=2,3,4$). \end{remark}


\begin{remark}\label{rem:consequences}
\begin{enumerate}
\item By Theorem \ref{TheoremB}, no non-circular planar Schiffer
domain satisfying \textup{(H1)} has
$\lambda|\Omega| \le \pi j_
{1,1}^2 \approx 46.1$.
\item If $u$ has a single critical point $x_0$ in $\Omega$, then
$\Omega$ is a disc. Indeed, since $\int_\Omega u = 0$ and
$u = c \neq 0$ on $\partial\Omega$, the function $u$ attains an
interior extremum, which must be $x
0$. Hence $u$ has no saddle
point, and $\Omega$ is a disc by the theorem of
Willms--Gladwell~\cite{WG94}. This improves the
convexity statement of Proposition~\ref{prop:cp}.
\end{enumerate}
\end{remark}
\section{Infinitesimal rigidity of the disc}\label{sec:rigidity}

\begin{theorem}\label{thm:A1}
Let $(\Om_\eps)$ be a family of Schiffer domains, $C^1$ in $\eps$ in $C^{2,\alpha}$, with $\Om_0$ the disc of
radius $R$ and $kR=j_{1,n}$. Then the first-order normal displacement is a combination of $1$, $\cos\theta$
and $\sin\theta$ (dilation and translations).
\end{theorem}

\begin{proof}
Write $u_\eps=u_0+\eps v$ and $\bOm_\eps=\{r=R(1+\eps h(\theta))\}$. Since $u_0'(R)=0$ and
$u_0''(R)=-k^2$, the conditions $u_\eps=1$ and $\dn u_\eps=0$ give, to first order, $v(R,\theta)=0$ and
$\partial_rv(R,\theta)=Rk^2h(\theta)$. For $h=\cos m\theta$ with $m\ge2$, we have $v=AJ_m(kr)\cos m\theta$, so
$AJ_m(kR)=0$. By Siegel's theorem~\cite{Si29} (Bourget's hypothesis), $J_1$ and $J_m$ have no common
positive zero; thus $J_m(j_{1,n})\neq0$, $A=0$ and $\partial_rv=0\neq Rk^2$, a contradiction. The mode $m=1$
corresponds to translations and $m=0$ to dilation.
\end{proof}
\begin{remark}
The above theorem is mainly in Canuto, \cite{Ca14}. Note that variations in $k$ and $c$
only affect the radial mode. We invite the reader to see also the references \cite{Ko93} and \cite{Ag93} on this subject. The authors established local rigidity in the neighborhood of the ball.

\end{remark}

This is consistent with~\cite{Ca14}: the counterexample of~\cite{CS26} does not branch off the Schiffer disc,
but off Wheeler's constant-flux family~\cite{Wh25}.

\begin{proposition}[Analogue for~\eqref{eq:berenstein}]\label{prop:A1B}
At the unit disc with $k=j_{0,n}$, a deformation $\cos m\theta$, $m\ge2$, is admissible to first order
for~\eqref{eq:berenstein} if and only if $kJ_m'(k)+J_m(k)=0$. This condition never holds for $m\geq 2, \forall n.$
\end{proposition}
\begin{proof}
With $u_0'(1)=1$ and $u_0''(1)=-1$, the linearized conditions are $v(1)=-h$ and $\partial_rv(1)=h$, hence
$AJ_m(k)=-1$ and $AkJ_m'(k)=1$. \\
By Lommel polynomials, $xJ_m^{'}+J_m= xJ_{m-1} - (m-1)J_m= S_m (x)J_1 (x)+ T_m (x)J_0 (x),$ where $S_m$ and $T_m$ are polynomials in $x$ and $1/x$ with rational coefficients.\\
For $x= j_{0, n}$ we have $J_0= 0$ and $J_1\neq 0:$ the condition is equivalent to $S_m (j_{0, n})= 0.$\\
$S_m$ is not null identically: its coefficient  in $x^{-(m-1)}$ equal $-(m-1) 2^{m-1} (m-1)!$.\\
$S_m (j_{0, n})= 0$ would make $j_{0, n}$ algebraic, but Siegel ( in 1929) showed that the roots of $J_0$  are transcendental.
\end{proof}

The disc is therefore infinitesimally rigid for problem $(2)$ as well. 

\begin{remark}[A curvature rigidity programme]
A quantitative version of Theorem~\ref{thm:A1}, namely an explicit $C^{2,\alpha}$ neighbourhood of size
$\eps(K)$ in which the disc is the only Schiffer domain with $kR^*\le K$, combined with quantitative
Alexandrov theorems~\cite{MP20,CM17}, would give: \emph{a planar Schiffer domain whose curvature oscillates
little is a disc.} The missing ingredient is a tame remainder estimate for the map
``shape $\mapsto$ flux''. 
\end{remark}

\subsection*{Acknowledgements and declaration}
The author would like to thank  colleagues from his research laboratory who encouraged him to interact with the AI Claude, Anthropic. 
The numerical computations and part of the writing were
carried out with the assistance of a language model (Claude, Anthropic). The author has
reviewed the whole text and takes full responsibility for it.



\begin{thebibliography}{xxxx}
\bibitem{Ag93} M. L. Agranovsky. On the stability of the spectrum in the Pompeiu problem. Journal of
Mathematical Analysis and Applications, 178(1):269- 279, 1993.
\bibitem{Av86}  Aviles P.; {\em Symmetry theorems related to Pompeiu's problem.} Amer. J. Math. 108 (5) (1986) 1023-1036. 

\bibitem{B} Berenstein C.A.; {\em An inverse spectral problem and its relation to the Pompeiu problem,} J. d'Analyse
Math. 37 (1980), 128-144.
\bibitem{B1} Berenstein C. A; {\em "The Pompeiu problem, What's new",} R. Deville (ed.) et al. (ed.) , Complex Analysis, Harmonic Analysis and Applications , Res. Notes Math. , 347 , Pitman (1996) pp. 1-11
\bibitem{BY}  Berenstein C. A., Yang P., {\em An overdetermined Neumann problem in the unit disc.} Adv. in Math., 44 (1982), 1-17.
\bibitem{BY87}  Berenstein C. A., Yang P., {\em An inverse  Neumann problem } J. Reine Angew. Math. 382 (1987), 1-21.
 \bibitem{BST73} Brown L,  Schreiber B.M. and  Taylor B.A.; {\em Spectral synthesis and the Pompeiu problem,}  Ann. Inst. Fourier (Grenoble) 23 (1973), 125-154.
\bibitem{Ca14} B.~Canuto, \emph{Stability results for the $N$-dimensional Schiffer conjecture via a perturbation method}, Calc. Var. PDE \textbf{50} (2014), 305-334.
\bibitem{CH99} Chatelain T. and Henrot A.; {\em Some results about Schiffer's conjectures } Inverse Problems 15 (1999) 647- 658.
\bibitem{CLDP26} G.~Cao-Labora, J.~de Dios Pont, \emph{Counterexamples to Schiffer's conjecture}, arXiv preprint, arXiv:2608.05114, August 2026
\bibitem{CLF25} G.~Cao-Labora, A.~J.~Fern\'andez, \emph{A contractible Schiffer counterexample on the half-sphere}, arXiv:2510.05732, 2025.
\bibitem{CM17} G.~Ciraolo, F.~Maggi, \emph{On the shape of compact hypersurfaces with almost constant mean curvature}, Comm. Pure Appl. Math. \textbf{70} (2017), 665-716.
\bibitem{CS26} M.~J.~Colbrook, G.~Stepaniants, \emph{A computer-assisted counterexample to the planar Pompeiu and Schiffer conjectures}, arXiv:2608.01579, 2026.
\bibitem{CSS26} M.~J.~Colbrook, S.~Sadeghi, G.~Stepaniants, \emph{A computer-assisted counterexample to the planar Berenstein conjecture}, arXiv:2608.08953, 2026.
 \bibitem{DeZo} Delfour M. C. and Zol\'esio J. P.; {\em  Shapes and Geometries: Metrics, Analysis, Differential
Calculus, and Optimization, } second edition, Society for Industrial and Applied Mathematics, 2011
\bibitem{DSWZ25} G.~Dai, Y.~Sun, J.~Wei, Y.~Zhang, \emph{On the Schiffer and Berenstein conjectures for centrally symmetric convex domains in the plane}, arXiv:2511.19819, 2025.

\bibitem{De12} Deng J.; {\em Some results on the Schiffer's conjecture in $\mathbb{R}^2$.} Journal of Differential Equations 253, no. 8 (2012): 2515-2526

\bibitem{Eb93} Ebenfelt P.; {\em  Some results on the Pompeiu problem,} Ann. Acad. Sci. Fennicae 18 (1993), 323-341. MR 95b:30061
\bibitem{FMW25} Fall M. M., Minlend I. A.  and Weth T.; {\em The Schiffer problem on the cylinder and on the 2-sphere}, J. Eur. Math. Soc., DOI 10.4171/JEMS/1612, (2025).
\bibitem{Gr00} A.~Grigelionis, \emph{On shape optimization and the Pompeiu problem}, arXiv:math/0012077, 2000.
\bibitem{Li07} Liu G.;{\em Symmetry theorems for the overdetermined eigenvalue problems.} J. Differential Equations 233 (2007), no. 2,
585- 600. 
\bibitem{Cvx26} J. Guo, \emph{ Convex counterexamples to the Schiffer and Pompeiu conjectures in dimensions two to
eighteen}, arXiv:2609.35419, v2 (2026)
\bibitem{Ko93} T. Kobayashi. Perturbation of domains in the Pompeiu problem. Communications in Analysis
and Geometry, 1(4):515-541, 1993
\bibitem{MM06} Michor, P.M., Mumford, D.; {\em Riemannian geometries on spaces of plane curves,} J. Eur. Math. Soc. 8(1), 1-48 (2006)
\bibitem{Mo26} Mondal S. {\em A short note on the Schiffer's conjecture for a class of centrally symmetric convex domains in $\mathbb{R}^2$,}  J. Anal. Math. 158 (2026), 401-412. 
\bibitem{MP20} R.~Magnanini, G.~Poggesi, \emph{Serrin's problem and Alexandrov's soap bubble theorem: enhanced stability via integral identities}, Indiana Univ. Math. J. \textbf{69} (2020), 1181-1205.
\bibitem{NSY20}  Nigam N., Siudeja B.,  Young B.; {\em A proof via finite elements for Schiffer's conjecture on a regular pentagon.} Found. Comput. Math. 20 (2020), no. 6, 1475-1504. 
\bibitem{Po29}  Pompeiu D.; {\em Sur une propri\'et\'e int\'egrale des fonctions de deux variables r\'eelles, Bull. Sci. Acad. Roy. Belgique 5 no 15 (1929) 265-269}
\bibitem{Se24} D.~Seck, \emph{On Pompeiu's--Schiffer's conjectures from shape optimization}, arXiv:2405.11700, 2024.
\bibitem{Se71} Serrin J.; {\em A Symmetry Problem in Potential Theory,} Arch. Rational Mech. Anal. 43, 304-318 (1971).
\bibitem{Si29} C.~L.~Siegel, \emph{\"Uber einige Anwendungen diophantischer Approximationen}, Abh. Preuss. Akad. Wiss. Phys.-Math. Kl. (1929), no.~1.
\bibitem{sozo}
Sokolowski, J. and  Zol\'esio, J. P.; {\em Introduction to Shape Optimization,}
Springer-Verlag, Paris, 1992.
\bibitem{Sp81} R.~P.~Sperb, \emph{Maximum Principles and Their Applications}, Academic Press, New York, 1981.
\bibitem{Ta76} G.~Talenti, \emph{Elliptic equations and rearrangements}, Ann. Scuola Norm. Sup. Pisa Cl. Sci. (4) \textbf{3} (1976), 697-718.
\bibitem{Ya82}  Yau S.T. ; {\em Seminars on Differential Geometry,} Ann. of Math. Stud., Princeton University Press,
1982.
\bibitem{Wh25} M.~H.~Wheeler, \emph{Non-symmetric solutions to an overdetermined problem for the Helmholtz equation in the plane}, arXiv:2509.00455, 2025.
\bibitem{WG94} N. B. Willms ,  G. M. L. Gladwell, \emph{ Saddle points and overdetermined problems for the Helmholtz equation} Z. Angew.
Math. Phys. \textbf{45}, pages 1-26 (1994)
 \bibitem{Wi76}  Williams S.A.; {\em A partial solution to the Pompeiu problem,} Math. Ann. 223, 1976,
183-190.
  \bibitem{Wi81}  Williams S.A.; {\em Analyticity of the boundary for Lipschitz domains without the Pompeiu property,}  Indiana Univ. Math. J. 30, 1981, 357-369 
.
\end{thebibliography}
\end{document}